\documentclass{amsart}
\usepackage{amssymb}

\def\squareforqed{\hbox{\rlap{$\sqcap$}$\sqcup$}}
\def\qed{\ifmmode\squareforqed\else{\unskip\nobreak\hfil
\penalty50\hskip1em\null\nobreak\hfil\squareforqed
\parfillskip=0pt\finalhyphendemerits=0\endgraf}\fi\medskip}

\newcommand{\udot}{{}^{\textstyle .}}

\newcommand{\PSL}{\mathrm{PSL}}
\newcommand{\PGL}{\mathrm{PGL}}

\newcommand{\Sz}{\mathrm{Sz}}

\newcommand{\Suz}{\mathrm{Suz}}
\newcommand{\PSU}{\mathrm{PSU}}

\newcommand{\Th}{\mathrm{Th}}
\newcommand{\MM}{\mathbb{M}}

\title{The uniqueness of $\PSU_3(8)$ in the Monster}
\date{Skeleton version begun 15th August 2015. First complete draft 18th July 2016.
This version 7th August 2016}
\author{Robert A. Wilson}

\address{School of Mathematical Sciences\\
Queen Mary University of London\\
London E1 4NS\\U.K.}
\email{R.A.Wilson@qmul.ac.uk}
\begin{document}
\maketitle

\begin{abstract}
As a  contribution to an eventual solution 
of the problem of the determination of the maximal subgroups of the Monster
we show that there is a unique conjugacy class of subgroups isomorphic to $\PSU_3(8)$.
The argument depends on some computations in various subgroups, but not
on computations in the Monster itself.
\end{abstract}

\section{Introduction}
The maximal subgroup problem for almost simple groups became a major focus for research
in group theory in the 1980s, and remains so today. In the case of the sporadic groups,
a systematic attack on the problem began earlier,
with Livingstone and his students in the 1960s.
The problem was solved in the 20th century
for $25$ of the $26$ sporadic
simple groups, and their automorphism groups, but one case, namely
the Fischer--Griess Monster group $\MM$, remains outstanding. 
A great deal of work on this case has already been done.
The maximal $p$-local
subgroups were classified in \cite{oddlocals,MSh,Meier}, and some theoretical work on
non-local subgroups was accomplished in \cite {Anatomy1,Anatomy2}.
Following successful computer constructions of the Monster \cite{3loccon,2loccon}
other techniques became available, and more progress was made
\cite{post,A5subs,S4subs,L241,L227,L213B,Sz8A}, 
including discovery of five previously unknown
maximal subgroups, isomorphic to 
\begin{itemize}
\item $\PSL_2(71)$, $\PSL_2(59)$,
$\PSL_2(41)$, $\PGL_2(29)$, $\PGL_2(19)$.
\end{itemize}

The cases left open by this published work are
possible maximal subgroups with socle isomorphic to one of the following simple groups:
\begin{itemize}
\item $\PSL_2(8)$, $\PSL_2(13)$, $\PSL_2(16)$, $\PSU_3(4)$, $\PSU_3(8)$.
\end{itemize}
Of these, $\PSL_2(8)$ and $\PSL_2(16)$ have been classified in unpublished work
of P. E. Holmes, although the results seem not to be publicly available. 
In this paper we deal with the case $\PSU_3(8)$.
Specifically, we show that, up to conjugacy, there is a unique
subgroup $\PSU_3(8)$ in the Monster. Its normalizer is the 
already known maximal subgroup $(A_5\times \PSU_3(8){:}3){:}2$.  
Notation 
follows \cite{Atlas,FSG}, where required
background information can be found.

\section{Existence}
Exactly one conjugacy class of subgroups of $\mathbb M$ isomorphic to
$\PSU_3(8)$ is contained in the known maximal subgroups. The normalizer of such
a group is $(A_5\times\PSU_3(8){:}3){:}2$, itself a maximal subgroup of $\mathbb M$.
For details, see \cite{Anatomy1}.

\section{Strategy for proving uniqueness}
The group $\PSU_3(8)$ can be generated from a group $\PSL_2(8)\times 3$ by 
extending $D_{18}\times 3$ to $(9\times 3){:}S_3$. Note that $9\times 3$
contains three cyclic subgroups of order $9$, which are permuted by the $S_3$. 
Similarly, there are three complements of order $3$, which are also permuted
by the $S_3$. Hence it is
sufficient to extend $9\times 3$ to $D_{18}\times 3$ normalizing one of the other
two cyclic subgroups of order $9$.

We note in particular that all cyclic groups of order $9$ in $\PSU_3(8)$
are conjugate, and hence we need only consider subgroups $\PSL_2(8)\times 3$
in which
the diagonal elements of order $9$ are
conjugate in the Monster
to the elements of order $9$ inside $\PSL_2(8)$. We shall show that there is
only one class of $\PSL_2(8)\times 3$ in the Monter
that satisfies this condition.
Moreover, the cyclic group of order $9$ extends to a unique
 $D_{18}$ in $\PSU_3(8)$. 
Hence the $D_{18}\times 3$ we wish to construct is conjugate
in the Monster to the one inside $\PSL_2(8)\times 3$.


\section{The subgroup $3\times\PSL_2(8)$}
Since $\PSL_2(8)$ contains elements of order $9$, the elements of order $3$
fuse to $\mathbb M$-class $3B$. Since it contains a pure $2^3$, the involutions are
in $\mathbb M$-class $2B$.
In \cite{Anatomy1} Norton accounts for many of the structure constants of type
$(2,3,7)$ in the Monster. In particular he shows that there is no $3\times\PSL_2(8)$ 
in which the $\PSL_2(8)$ is of type
$(2B,3B,7B)$. He also shows that there are three classes of $\PSL_2(8)$
of type $(2B,3B,7A)$, just two of which centralize elements of order $3$.
The respective normalizers are:
\begin{enumerate}
\item $\PSL_2(8){:}3\times 3S_6$.
Here the central $3$ in $3A_6$ is in Monster class $3A$, as are the $3$-cycles.
The elements mapping to fixed-point-free $3$-elements in $A_6$ are in Monster class $3C$.
\item $\PSL_2(8)\times 2\times S_4$. Here, the
elements of order $3$ in $S_4$ are in Monster class $3A$.
\end{enumerate}
Hence there are exactly four classes of $\PSL_2(8)\times 3$ in the Monster.
%

\section{Fusion of elements of order $9$}
Consider first the case where the central elements of $3\times \PSL_2(8)$ are in
class $3C$ in the Monster.
We restrict the character of degree $196883$ to $S_3\times \Th$. Using the character
values on $3C$ and $2A$, we obtain a decomposition as
$$2\otimes 65628 + 1^+\otimes 34999 + 1^-\otimes 30628$$
where the first factor denotes the representation of $S_3$. The values on classes $2A$
and $7A$ of $\Th$ are easily computed:
$$\begin{array}{rrr}
1A & 2A & 7A\cr
34999 & 183 & 13\cr
30628 & -92 & 3\cr
65628 & 92 & 17
\end{array}$$
from which it is easy to see that the decomposition into irreducibles of $\Th$ is given by
\begin{eqnarray*}
34999 &=& 30875+4123+1\\
30628&=& 30628\\
65628&=& 61256+4123+248+1
\end{eqnarray*}


It then follows that the values of the character of degree $196883$ on elements of
$\Th$-class $9A$, $9B$, and $9C$ are respectively $-1$, $-1$, $26$, while the
values on the corresponding diagonal elements in $3\times\Th$ are $26$, $26$, and
$-1$ respectively. In other words, the diagonal elements are always in a 
different conjugacy class from the elements in $\Th$. Hence this case is
eliminated. (In fact, in this case the $\PSL_2(8)$ contains elements of $\Th$ class $9C$,
that is, Monster class $9A$.)

The remaining three classes of $\PSL_2(8)\times 3$, namely the ones with a central
$3A$-element, are contained in the double cover of
the Baby Monster. The work in \cite{maxB} then shows that in these cases the 
elements of order $9$ in $\PSL_2(8)$ are in Baby Monster class $9B$, so
Monster class $9A$. Moreover, in two of the three cases, the diagonal elements
of order $9$ are in Baby Monster class $9A$, so Monster class $9B$. But in $\PSU_3(8)$
these two classes of elements of order $9$ are fused. Hence these cases cannot
extend to $\PSU_3(8)$ in the Monster.

The remaining case therefore is
a $3A$-type, with normalizer $\PSL_2(8){:}3\times 3S_6$.
We know there exists such a subgroup
$\PSU_3(8)$ in the Monster, so all
elements of order $9$ fuse to $9A$.

\section{The centralizer of a $9A$ element}
From \cite{oddlocals}, the centralizer of a $9A$-element in the Monster has shape
$[3^7].\PSU_4(2)$. Looking more closely, we see that the structure is the central product of
the cyclic group of order $9$ with a group of shape $3^6\udot \Omega_5(3)$, in which the
action of $\Omega_5(3)$ on $3^6$ is uniserial, with a trivial submodule and a natural
module as quotient.
Moreover, since
this group contains $9 \times 3\udot S_6$, 
the extension is non-split, in the sense that $C(9)/9\cong 3^5\udot \Omega_5(3)$.

These facts can be checked computationally, using the construction of the subgroup
$3^{1+12}\udot2\udot\Suz{:}2$ described 
in \cite{3loccon}. But in fact the proof below does not
depend on any of the subtleties, so the sceptical reader can ignore them.

\section{The centralizer of $9\times 3$}
Centralizing the additional element of order $3$ reduces the group from
$9\circ3^6\udot \Omega_5(3)$ to $9\circ3^6\udot A_6$. The structure of the latter group
is very subtle, and in particular it contains several conjugacy classes of $3\udot A_6$,
and it is not obvious which one centralizes $\PSL_2(8)$.

In any case, the group of elements which either centralize $9\times 3$ or extend it to
$D_{18}\times 3$ is of shape $(9\circ 3^6)\udot(A_6\times 2)=
(9\times 3).3^4.(2\times A_6)$. We must adjoin an
involution in the conjugacy class which maps to the central involution in the
quotient $A_6\times 2$. But this conjugacy class contains only $3^6=729$ elements,
while the group of symmetries is $9\times 3A_6$, of order $9720$. Hence every
group generated in the prescribed fashion has non-trivial centralizer in the Monster.
Indeed, this counting argument implies that such a centralizer has order at least 
$9720/729=13\frac13$.

\section{Proof of the theorem}
The centralizer of an element of order $19$ is $19\times A_5$, containing elements of
classes $2A$, $3C$ and $5A$. The only subgroup of $A_5$ with order at least $14$
is $A_5$ itself. Hence every $\PSU_3(8)$ in the Monster has centralizer conjugate to this
$A_5$.

As a corollary we obtain new proofs of the uniqueness of $\PSU_3(8)$ as a subgroup of
the Baby Monster, the Thompson group, and the Harada--Norton group.


\end{document}